\documentclass[12pt]{article}
\usepackage{amsmath}
\usepackage{amssymb}
\usepackage{vatola}

\def\q{\quad}

\def\qtq#1{\q\t{#1}\q}
\def\mod#1{\ (\text{\rm mod}\ #1)}
\def\t{\text}
\def\f{\frac}
\def\e{\equiv}
\def\b{\binom}

\def\sls#1#2{(\f{#1}{#2})}
 
\def\Ls#1#2{\Big(\f{#1}{#2}\Big)}

\let \pro=\proclaim
\let \endpro=\endproclaim

\begin{document}
%\leftline{This is the revision for possible publication in Integral
%Transforms and Special Functions.} \leftline{(February 16, 2015)}
\par\q\par\q
 \centerline {\bf
Congruences for Domb and Almkvist-Zudilin numbers}
$$\q$$
\centerline{Zhi-Hong Sun} $$\q$$ \centerline{School of Mathematical
Sciences, Huaiyin Normal University,} \centerline{Huaian, Jiangsu
223001, P.R. China} \centerline{Email: zhihongsun@yahoo.com}
\centerline{Homepage: http://www.hytc.edu.cn/xsjl/szh}
 \abstract{In this paper we prove some transformation formulae for
congruences modulo a prime and deduce some congruences for Domb
numbers and Almkvist-Zudilin numbers. We also pose some conjectures
on congruences modulo prime powers.
\par\q
 }
 \endabstract
\let\thefootnote\relax \footnotetext {The author is supported by
the Natural Science Foundation of China (grant No. 11371163).}

\section*{1. Introduction}
\par\q  Let $[x]$ be the greatest integer not exceeding $x$,
and let $\sls ap$ be the Legendre symbol. For a prime $p$ let $\Bbb
Z_p$ be the set of rational numbers whose denominator is not
divisible by $p$. For positive integers $a,b$ and $n$, if
$n=ax^2+by^2$ for some integers $x$ and $y$, we briefly write that
$n=ax^2+by^2$.
\par Let $(a)_0=1$ and $(a)_k=a(a+1)\cdots(a+k-1)$ for any positive integer $k$.
 Then $\f{(a)_k}{k!}
=(-1)^k\b{-a}k$. A formula of Bailey (see [GZ, (9) and (12)]) states
that
$$\align \sum_{k=0}^{\infty}\f{\sls 12_k^3}{k!^3}x^k&=\f
2{\sqrt{4-x}}\sum_{k=0}^{\infty}\f{\sls 12_k\sls 16_k\sls
56_k}{k!^3} \Ls{27x^2}{(4-x)^3}^k
\\&=\f
1{\sqrt{1-4x}}\sum_{k=0}^{\infty}\f{\sls 12_k\sls 16_k\sls
56_k}{k!^3} \Big(\f{27x}{(4x-1)^3}\Big)^k,
\endalign$$
where $|x|$ is sufficiently small. It is easily seen that
$$\f{\sls 12_k}{k!}=\f{\b{2k}k}{4^k}
\qtq{and} \f{\sls 12_k\sls 16_k\sls 56_k}{k!^3}
=\f{\b{2k}k\b{3k}k\b{6k}{3k}}{12^{3k}}.\tag 1.1$$ Thus, taking
$x=64/m$ $(|m|\ \t{is sufficiently large})$ in Bailey's
transformation and applying (1.1) we get
$$\aligned\sum_{k=0}^{\infty}\f{\b{2k}k^3}{m^k}
&=\sqrt{\f m{m-16}}\sum_{k=0}^{\infty}\b{2k}k\b{3k}k\b{6k}{3k} \Ls
m{(m-16)^3}^k
\\&=\sqrt{\f m{m-256}}\sum_{k=0}^{\infty}\b{2k}k\b{3k}k\b{6k}{3k} \Ls
{m^2}{(256-m)^3}^k.\endaligned\tag 1.2$$ Let $p$ be an odd prime and
$m\in\Bbb Z_p$ with $m\not\e 0\mod p$. In [Su1], Z.W. Sun
conjectured many congruences modulo $p^2$ for the sums
$$\sum_{k=0}^{p-1}\f{\b{2k}k^3}{m^k},\
\sum_{k=0}^{p-1}\f{\b{2k}k^2\b{3k}k}{m^k},\
\sum_{k=0}^{p-1}\f{\b{2k}k^2\b{4k}{2k}}{m^k} \qtq{and}
\sum_{k=0}^{p-1}\f{\b{2k}k\b{3k}k\b{6k}{3k}}{m^k}.$$ In [S2-S4] the
author confirmed some of his conjectures. In Section 2, using some
results in [S3, S4] we prove the following p-analogue of (1.2):
$$\aligned\sum_{k=0}^{(p-1)/2}\f{\b{2k}k^3}{m^k}
&\e \Ls {m(m-16)}p\sum_{k=0}^{[p/6]}\b{2k}k\b{3k}k\b{6k}{3k} \Ls
m{(m-16)^3}^k
\\&\e \Ls {m(m-256)}p\sum_{k=0}^{[p/6]}\b{2k}k\b{3k}k\b{6k}{3k} \Ls
{m^2}{(256-m)^3}^k \mod p.\endaligned$$  We also obtain similar
congruences for
$\sum_{k=0}^{[p/3]}\b{2k}k^2\b{3k}k\sls{1-t^2}{108}^k$ and
$\sum_{k=0}^{[p/4]}\b{2k}k^2\b{4k}{2k}\sls{1-t^2}{256}^k$.
\par For any nonnegative integer n let
$$\aligned &D_n=\sum_{k=0}^n\b{2k}k\b{2n-2k}{n-k}\b nk^2,
\q a_n=\sum_{k=0}^n\b nk^2\b{2k}k,\\&
b_n=\sum_{k=0}^{[n/3]}\b{2k}k\b{3k}k\b n{3k}\b{n+k}k(-3)^{n-3k}.
\endaligned\tag 1.3$$
The first few values of $\{D_n\},\{a_n\}$ and $\{b_n\}$ are as
below:
$$\align &D_0=1,\ D_1=4,\ D_2=28,\ D_3=256,\ D_4=2716,\ D_5=31504,
\ D_6=387136,
\\&a_0=1,\ a_1=3,\ a_2=15,\ a_3=93,\ a_4=639,\ a_5=4653,\ a_6=35169,
\\&b_0=1,\ b_1=-3,\ b_2=9,\ b_3=-3,\ b_4=-279,\ b_5=2997,\
b_6=-19431.\endalign$$ The numbers $\{D_n\}$ are called Domb numbers
since Domb introduced it in 1960, and the numbers $\{b_n\}$ are
called Almkvist-Zudilin numbers. See [CCL], [CV], [CZ], [D], [Z],
and A002895, A002893 and A125143 in A. Sloane's ``The on-line
encyclopedia of integer sequences". Such sequences appear as
coefficients in various series for $1/\pi$. For example, from [CCL]
and [CV] we know that
$$\sum_{n=0}^{\infty}\f{5n+1}{64^n}D_n=\f 8{\sqrt 3\pi}\qtq{and}
\sum_{n=0}^{\infty}\f{4n+1}{81^n}b_n=\f{3\sqrt 3}{2\pi}.$$ In [R],
by using very advanced and complicated method Rogers showed that
$$\sum_{n=0}^{\infty}D_nu^n
=\f
1{1-4u}\sum_{k=0}^{\infty}\b{2k}k^2\b{3k}k\Ls{u^2}{(1-4u)^3}^k\tag
1.4$$ and $$\sum_{k=0}^{\infty}\b{2k}k\Ls u{9(1+u)^2}^ka_k
=\f{1+u}{1+3u}\sum_{k=0}^{\infty}\b{2k}k^2\b{4k}{2k}\Ls
u{9(1+3u)^4}^k,\tag 1.5$$
 where $|u|$ is sufficiently small.
\par
Let $p$ be an odd prime and $u\in\Bbb Z_p$. In Sections 3 and 4 we
prove that
$$\sum_{n=0}^{p-1}D_nu^n
\e\sum_{k=0}^{p-1}\b{2k}k^2\b{3k}k\Ls{u^2}{(1-4u)^3}^k\mod p\q\t{for
$u\not\e \f 14\mod p$},$$ and for $u\not\e -\f 19,-\f 1{27}\mod p$,
$$\sum_{k=0}^{p-1}\b{2k}ka_k\Ls u{(1+9u)^2}^k\e
\sum_{n=0}^{p-1}b_nu^n \e\sum_{k=0}^{p-1}\b{2k}k^2\b{4k}{2k}\Ls
u{(1+27u)^4}^k\mod p,$$ which are $p$-analogues of (1.4) and (1.5).
 As an application we prove congruences for $\sum_{n=0}^{p-1}
 \f{D_n}{m^n}$ modulo $p$ for $m=1,-2,4,8,-8,16,-32,64$, which were conjectured by the author's brother Z.W. Sun
in [Su3]. For instance, if $p\e 1,4\mod{15}$ is a prime and so
$p=x^2+15y^2$, then $\sum_{n=0}^{p-1}D_n\e 4x^2\mod p.$ In Sections
4 and 5 we determine $\sum_{k=0}^{p-1}\b{2k}k\f{a_k}{m^k}\mod p$ for
$m=-12,36,100$, and $\sum_{n=0}^{p-1}\f{b_n}{m^n}\mod p$ for
$m=1,-3,9,-9,$ $-27,81$.
 We also determine
 $\sum_{k=0}^{p-1}\frac{\b{2k}k^2\b{3k}k}{1458^k}, \
 \sum_{k=0}^{p-1}\frac{\b{2k}k^2\b{3k}k}{15^{3k}}$ and
 $\sum_{k=0}^{p-1}\frac{\b{2k}k^2\b{4k}{2k}}{28^{4k}}\mod p$
 and so partially confirm three conjectures in [S1] and [Su1].

 \par In Section 6, we pose some conjectures on congruences
 modulo prime powers.

 \section*{2. Transformation formulas involving $\b{2k}k\b{3k}k\b{6k}{3k}$}
\par\q Let $p$ be an odd prime and $k\in\{0,1,\ldots,p-1\}$. It is easily
seen that (see [S2-S3])
$$\align &\b{2k}k=\f{(2k)!}{k!^2}\e 0\mod p\qtq{for}k>\f p2,
\\&\b{2k}k\b{3k}k=\f{(3k)!}{k!^3}\e 0\mod p\qtq{for}k>\f p3,
\\&\b{2k}k\b{4k}{2k}=\f{(4k)!}{(2k)!\cdot k!^2}\e 0\mod p
\qtq{for}k>\f p4,
\\&\b{3k}k\b{6k}{3k}=\f{(6k)!}{(3k)!(2k)!k!}\e 0\mod p
\qtq{for}k>\f p6.\endalign$$

\par Let $\{P_n(x)\}$ be the Legendre polynomials given by
 $$P_n(x)=\f
1{2^n}\sum_{k=0}^{[n/2]}\b nk(-1)^k\b{2n-2k}nx^{n-2k} =\f 1{2^n\cdot
n!}\cdot\f{d^n}{dx^n}(x^2-1)^n.$$ Then clearly
$P_n(-x)=(-1)^nP_n(x)$.
 In [S2, Theorems 3.1 and 4.1] the
author showed that for any prime $p>3$ and $t\in \Bbb Z_p$,
$$P_{[\f  p3]}(t) \e -\Ls
p3\sum_{x=0}^{p-1}\Ls{x^3+3(4t-5)x+2(2t^2-14t+11)}p\mod p\tag 2.1$$
and
$$\sum_{k=0}^{[p/3]}\b{2k}k^2\b{3k}k\Ls{1-t^2}{108}^k\e P_{[\f
p3]}(t)^2\mod p.\tag 2.2$$ In [S3, Theorems 2.1 and 4.2] the author
showed that for any prime $p>3$ and $t\in \Bbb Z_p$,
$$P_{[\f p4]}(t)\e -\Ls{6}p\sum_{x=0}^{p-1}\Ls{x^3-\f{3(3t+5)}2x+9t+7}p\mod
 p\tag 2.3$$
and
$$\sum_{k=0}^{[p/4]}\b{2k}k^2\b{4k}{2k}\Ls{1-t^2}{256}^k\e P_{[\f
p4]}(t)^2\mod p.\tag 2.4$$ In [S4, Theorem 4.2], the author showed
that for any prime $p>3$ and $m,n\in \Bbb Z_p$ with $m\not\e 0\mod
p$,
$$\Big(\sum_{x=0}^{p-1}\Big(\frac{x^3+mx+n}p\Big)\Big)^2
\e \Big(\frac{-3m}p\Big)
\sum_{k=0}^{[p/6]}\binom{2k}k\binom{3k}k\binom{6k}{3k}
\Big(\frac{4m^3+27n^2}{12^3\cdot 4m^3}\Big)^k\mod p.\tag 2.5$$
\pro{Theorem 2.1} For any prime $p>3$ and $m\in\Bbb Z_p$ with
$m\not\e 0,16,64,256\mod p$ we have
$$\aligned\sum_{k=0}^{(p-1)/2}\f{\b{2k}k^3}{m^k}
&\e \Ls {m(m-16)}p\sum_{k=0}^{[p/6]}\b{2k}k\b{3k}k\b{6k}{3k} \Ls
m{(m-16)^3}^k
\\&\e \Ls {m(m-256)}p\sum_{k=0}^{[p/6]}\b{2k}k\b{3k}k\b{6k}{3k} \Ls
{m^2}{(256-m)^3}^k \mod p.\endaligned$$\endpro
 Proof. By [S3, Theorem 3.2],
$$\sum_{k=0}^{(p-1)/2}\f{\b{2k}k^3}{m^k} \e \Ls{m(m-64)}pP_{[\f
p4]}\Ls{m+64}{m-64}^2\mod p.$$ Thus, applying (2.3) and (2.5) we
deduce that
$$\align \sum_{k=0}^{(p-1)/2}\f{\b{2k}k^3}{m^k}
&\e\Ls {m(m-64)}p\Big(\sum_{x=0}^{p-1}\Ls{x^3-\f
32(3\cdot\f{m+64}{m-64}+5)x+9\cdot\f{m+64}{m-64}+7}p\Big)^2
\\&\e \Ls {m(m-16)}p\sum_{k=0}^{[p/6]}\b{2k}k\b{3k}k\b{6k}{3k}
\Ls m{(m-16)^3}^k \mod p.\endalign$$ As $P_n(-x)=(-1)^nP_n(x)$, we
also have
$$\align \sum_{k=0}^{(p-1)/2}\f{\b{2k}k^3}{m^k}
&\e  \Ls{m(m-64)}pP_{[\f p4]}\Big(-\f{m+64}{m-64}\Big)^2
\\&\e\Ls {m(m-64)}p\Big(\sum_{x=0}^{p-1}\Ls{x^3-\f
32(-3\cdot\f{m+64}{m-64}+5)x-9\cdot\f{m+64}{m-64}+7}p\Big)^2
\\&\e \Ls {m(m-256)}p\sum_{k=0}^{[p/6]}\b{2k}k\b{3k}k\b{6k}{3k}
\Ls {m^2}{(256-m)^3}^k \mod p.\endalign$$
This proves the theorem.

\pro{Theorem 2.2} For any prime $p>3$ and $t\in\Bbb Z_p$ with
$4t\not\e \pm 5\mod p$ we have
$$\align \sum_{k=0}^{[p/3]}\b{2k}k^2\b{3k}k\Ls{1-t^2}{108}^k
&\e \Ls{5-4t}p\sum_{k=0}^{[p/6]}\b{2k}k\b{3k}k\b{6k}{3k}
\Ls{(t-1)(t+1)^3}{432(4t-5)^3}^k
\\&\e
\Ls{5+4t}p\sum_{k=0}^{[p/6]}\b{2k}k\b{3k}k\b{6k}{3k}
\Ls{(t+1)(1-t)^3}{432(4t+5)^3}^k \mod p.\endalign$$
\endpro
Proof. By (2.1), (2.2) and (2.5),
$$\align&\sum_{k=0}^{[p/3]}\b{2k}k^2\b{3k}k\Ls{1-t^2}{108}^k
\\&\e P_{[\f p3]}(t)^2\e \Big(\sum_{x=0}^{p-1}\Ls{x^3+3(4t-5)x+2(2t^2-14t+11)}p
\Big)^2
\\&\e \Ls{5-4t}p\sum_{k=0}^{[p/6]}\b{2k}k\b{3k}k\b{6k}{3k}
\Ls{(t-1)(t+1)^3}{432(4t-5)^3}^k\mod p.\endalign$$ Substituting $t$
with $-t$ in the above congruence we obtain the remaining result.
\par\q\par{\bf Remark 2.1} Taking $t=2z-1$ in Theorem 2.2 we see
that for any prime $p>3$ and $z\not\e\f 98\mod p$,
$$\sum_{k=0}^{[p/3]}\b{2k}k^2\b{3k}k\Ls{z(1-z)}{27}^k
\e \Ls{9-8z}p\sum_{k=0}^{[p/6]}\b{2k}k\b{3k}k\b{6k}{3k}
\Ls{z^3(1-z)}{27(9-8z)^3}^k\mod p.$$ This can be viewed as the
p-analogue of the Kummer-Coursat transformation ([GZ, (20)]):
$$\sum_{k=0}^{\infty}\b{2k}k^2\b{3k}k\Ls{z(1-z)}{27}^k
=\f 3{\sqrt{9-8z}}\sum_{k=0}^{\infty}\b{2k}k\b{3k}k\b{6k}{3k}
\Ls{z^3(1-z)}{27(9-8z)^3}^k.$$

 \pro{Theorem
2.3} For any prime $p>3$ and $t\in\Bbb Z_p$ with $3t\not\e \pm 5\mod
p$, we
 have
$$\align\sum_{k=0}^{[p/4]}\b{2k}k^2\b{4k}{2k}\Ls{1-t^2}{256}^k
&\e
\Ls{10+6t}p\sum_{k=0}^{[p/6]}\b{2k}k\b{3k}k\b{6k}{3k}\Ls{(t-1)^2(t+1)}
{64(3t+5)^3}^k
\\&\e
\Ls{10-6t}p\sum_{k=0}^{[p/6]}\b{2k}k\b{3k}k\b{6k}{3k}\Ls{(t+1)^2(t-1)}
{64(3t-5)^3}^k \mod p.\endalign$$ \endpro
 Proof. By (2.3), (2.4) and
(2.5),
$$\align\sum_{k=0}^{[p/4]}\b{2k}k^2\b{4k}{2k}\Ls{1-t^2}{256}^k
&\e P_{[\f p4]}(t)^2\e
\Big(\sum_{x=0}^{p-1}\Ls{x^3-\f{3(3t+5)}2x+9t+7}p \Big)^2
\\&\e \Ls{10+6t}p\sum_{k=0}^{[p/6]}\b{2k}k\b{3k}k\b{6k}{3k}
\Ls{(t-1)^2(t+1)}{64(3t+5)^3}^k\mod p.\endalign$$ Substituting $t$
with $-t$ in the above congruence we obtain the remaining result.

 \section*{3. Congruences  involving $\{D_n\}$}
 \pro{Lemma 3.1} Let $n$ be a nonnegative integer. Then
$$\sum_{k=0}^{[n/2]}\b{2k}k^2\b{3k}k\b{n+k}{3k}4^{n-2k}
=\sum_{k=0}^n\b{2k}k\b{2n-2k}{n-k}\b
nk^2.$$
\endpro
Proof. Let $S_1(n)$ and $S_2(n)$ denote the left side and the right
side of the identity, respectively. Using Maple and
 the D. Zeilberger¡¯s MAPLE programme EKHAD (Zeilberger algorithm)
  we find that for $i=1,2$,
$$(m+2)^3S_i(m+2)-2(2m+3)(5m^2+15m+12)S_i(m+1)
+ 64(m+1)^3S_i(m)=0\q (m=0,1,\ldots).$$ Since $S_1(0)=1=S_2(0)$ and
$S_1(1)=4=S_2(1)$, we deduce that $S_1(n)=S_2(n)$ for all
$n=0,1,2,\ldots$. This completes the proof.
\par\q\par
{\bf Proof of (1.4)}: By Lemma 3.1,
$$\align &\f 1{1-4u}\sum_{k=0}^{\infty}\b{2k}k^2\b{3k}k\Ls{u^2}{(1-4u)^3}^k
=\sum_{k=0}^{\infty}\b{2k}k^2\b{3k}ku^{2k}(1-4u)^{-3k-1}
\\&=\sum_{k=0}^{\infty}\b{2k}k^2\b{3k}ku^{2k}\sum_{r=0}^{\infty}
\b{-3k-1} r(-4u)^r
\\&=\sum_{n=0}^{\infty}u^n\sum_{k=0}^{[n/2]}\b{2k}k^2\b{3k}k
\b{-3k-1}{n-2k}(-4)^{n-2k}
\\&=\sum_{n=0}^{\infty}u^n\sum_{k=0}^{[n/2]}\b{2k}k^2
\b{3k}k\b{n+k}{3k}4^{n-2k} =\sum_{n=0}^{\infty}D_nu^n.\endalign$$

 \pro{Theorem 3.1} Let $p$ be an odd
prime and $u\in\Bbb Z_p$ with $u\not\e \f 14\mod p$. Then
$$\sum_{n=0}^{p-1}D_nu^n
\e\sum_{k=0}^{p-1}\b{2k}k^2\b{3k}k\Ls{u^2}{(1-4u)^3}^k\mod p.$$
\endpro
Proof. As $p\mid \b{2k}k\b{3k}k$ for $\f p3<k<p$, using Fermat's
little theorem and Lemma 3.1 we see that
$$\align &\sum_{k=0}^{p-1}\b{2k}k^2\b{3k}k\Ls{u^2}{(1-4u)^3}^k
\e \sum_{k=0}^{[p/3]}\b{2k}k^2\b{3k}ku^{2k}(1-4u)^{p-1-3k}
\\&=\sum_{k=0}^{[p/3]}\b{2k}k^2\b{3k}ku^{2k}\sum_{r=0}^{p-1-3k}\b{p-1-3k}r(-4u)^r
\\&\e\sum_{n=0}^{p-1}u^n\sum_{k=0}^{[n/2]}\b{2k}k^2\b{3k}k\b{p-1-3k}{n-2k}(-4)^{n-2k}
\\&\e \sum_{n=0}^{p-1}u^n\sum_{k=0}^{[n/2]}\b{2k}k^2\b{3k}k\b{-1-3k}{n-2k}(-4)^{n-2k}
\\&=\sum_{n=0}^{p-1}u^n\sum_{k=0}^{[n/2]}\b{2k}k^2\b{3k}k\b{n+k}{3k}4^{n-2k}
=\sum_{n=0}^{p-1}D_nu^n\mod p.\endalign$$ Thus the theorem is
proved.

 \pro{Theorem 3.2} Let $p$ be a
prime such that $p\e 1,4\mod 5$. Then
$$\align \sum_{n=0}^{p-1}D_n
\e\cases 4x^2\mod p&\t{if $p\e 1,4\mod {15}$ and so $p=x^2+15y^2$,}
\\0\mod p&\t{if $p\e 11,14\mod{15}$.}
\endcases\endalign$$
\endpro
Proof. Taking $u=1$ in Theorem 3.1 and then applying [S2, Theorem
4.6] we obtain the result.
\par{\bf Remark 3.1} In [Su3], Z.W. Sun conjectured that
for any prime $p>5$,
$$\align &\Ls p3\sum_{n=0}^{p-1}\sum_{k=0}^n\b nk^4
\e\sum_{n=0}^{p-1}D_n\e \sum_{n=0}^{p-1}\f{D_n}{64^n} \\&\e\cases
4x^2-2p\mod{p^2}&\t{if $p=x^2+15y^2\e 1,4\mod{15}$,}
\\2p-12x^2\mod{p^2}&\t{if $p=3x^2+5y^2\e 2,8\mod{15}$,}
\\0\mod{p^2}&\t{if $p\e 7,11,13,14\mod{15}$.}
\endcases\endalign$$

\pro{Theorem 3.3} Let $p$ be a prime such that $p\e 1,7,17,23\mod
{24}$. Then
$$\sum_{n=0}^{p-1}\f {D_n}{(-8)^n}
\e\cases 4x^2\mod p&\t{if $p\e 1,7\mod {24}$ and so $p=x^2+6y^2$,}
\\0\mod p&\t{if $p\e 17,23\mod{24}$.}
\endcases$$
\endpro
Proof. Taking $u=-\f 18$ in Theorem 3.1 and then applying [S2,
Theorem 4.5] we obtain the result.

\par{\bf Remark 3.2} In [Su3], Z.W. Sun conjectured that
for any prime $p>3$,
$$\sum_{n=0}^{p-1}\f{D_n}{(-8)^n}\e
\cases 4x^2-2p\mod{p^2}&\t{if $p=x^2+6y^2\e 1,7\mod{24}$,}
\\8x^2-2p\mod{p^2}&\t{if $p=2x^2+3y^2\e 5,11\mod{24}$,}
\\0\mod{p^2}&\t{if $p\e 13,17,19,23\mod{24}$.}
\endcases$$
\pro{Theorem 3.4} Let $p$ be an odd prime. Then
$$\sum_{n=0}^{p-1}\f {D_n}{8^n}
\e\cases 4x^2\mod p&\t{if $p\e 1,3\mod {8}$ and so $p=x^2+2y^2$,}
\\0\mod p&\t{if $p\e 5,7\mod 8$.}
\endcases$$
\endpro
Proof. Taking $u=\f 18$ in Theorem 3.1 and then applying [S2,
Theorem 4.3] we obtain the result.
\par{\bf Remark 3.3} In [Su3], Z.W. Sun conjectured that
for any odd prime $p$,
$$\sum_{n=0}^{p-1}\f{D_n}{8^n}\e
\cases 4x^2-2p\mod{p^2}&\t{if $p=x^2+2y^2\e 1,3\mod 8$,}
\\0\mod{p^2}&\t{if $p\e 5,7\mod 8$.}
\endcases$$

\pro{Lemma 3.2 ([CZ, Corollary 3.4]} Let $n$ be a nonnegative
integer. Then
$$D_n=\sum_{k=0}^n(-1)^k16^{n-k}\b nk\b{2k}k^2\b{n+2k}n.$$
\endpro
\par Lemma 3.2 can also be proved by using Maple and
 the D. Zeilberger¡¯s MAPLE programme EKHAD.

\pro{Theorem 3.5} Let $p$ be an odd prime and $u\in\Bbb Z_p$ with
$u\not\e \f 1{16}\mod p$. Then
$$\sum_{n=0}^{p-1}D_nu^n
\e\sum_{k=0}^{p-1}\b{2k}k^2\b{3k}k\Ls{-u}{(1-16u)^3}^k\mod p.$$
\endpro
Proof. As $p\mid \b{2k}k\b{3k}k$ for $\f p3<k<p$, using Fermat's
little theorem and Lemma 3.2 we see that
$$\align &\sum_{k=0}^{p-1}\b{2k}k^2\b{3k}k\Ls{-u}{(1-16u)^3}^k
\e \sum_{k=0}^{[p/3]}\b{2k}k^2\b{3k}k(-u)^k(1-16u)^{p-1-3k}
\\&=\sum_{k=0}^{[p/3]}\b{2k}k^2\b{3k}k(-u)^k
\sum_{r=0}^{p-1-3k}\b{p-1-3k}r(-16u)^r
\\&\e\sum_{n=0}^{p-1}u^n\sum_{k=0}^n\b{2k}k^2\b{3k}k(-1)^k
\b{p-1-3k}{n-k}(-16)^{n-k}
\\&\e \sum_{n=0}^{p-1}u^n\sum_{k=0}^n
\b{2k}k^2\b{3k}k(-1)^k\b{-1-3k}{n-k}(-16)^{n-k}
\\&=\sum_{n=0}^{p-1}u^n\sum_{k=0}^n
\b{2k}k^2\b{3k}k(-1)^k\b{n+2k}{n-k}16^{n-k}
\\&=\sum_{n=0}^{p-1}u^n\sum_{k=0}^n \b{2k}k^2\b nk
\b{n+2k}n(-1)^k16^{n-k} =\sum_{n=0}^{p-1}D_nu^n\mod p.\endalign$$
Thus the theorem is proved.
 \pro{Corollary 3.1} Let $p$ be an odd prime,
$u\in\Bbb Z_p$ and $u\not\e \f 14,\f 1{16}\mod p$. Then
$$\sum_{k=0}^{p-1}\b{2k}k^2\b{3k}k\Ls{-u}{(1-16u)^3}^k
\e\sum_{k=0}^{p-1}\b{2k}k^2\b{3k}k\Ls{u^2}{(1-4u)^3}^k\mod p.$$
\endpro
Proof. This is immediate from Theorems 3.1 and 3.5.

\par Corollary 3.1 is the $p$-analogue of the following formula in
[R, equation (3.6)]:
$$\f 1{1-16u}\sum_{k=0}^{\infty}\b{2k}k^2\b{3k}k\Ls{-u}{(1-16u)^3}^k
=\f
1{1-4u}\sum_{k=0}^{\infty}\b{2k}k^2\b{3k}k\Ls{u^2}{(1-4u)^3}^k.$$

\pro{Theorem 3.6} Let $p$ be a prime such that $p\e 1,4\mod 5$. Then
$$\align \sum_{n=0}^{p-1}\f{D_n}{64^n}
\e \sum_{k=0}^{p-1}\f{\b{2k}k^2\b{3k}k}{15^{3k}}
 \e\cases 4x^2\mod
p&\t{if $p\e 1,4\mod {15}$ and so $p=x^2+15y^2$,}
\\0\mod p&\t{if $p\e 11,14\mod{15}$.}
\endcases\endalign$$
\endpro
Proof. Taking $u=\f 1{64}$ in Theorem 3.5 and Corollary 3.1 we see
that
$$\sum_{n=0}^{p-1}\f{D_n}{64^n}\e
\sum_{k=0}^{p-1}\f{\b{2k}k^2\b{3k}k}{(-27)^k} \e
\sum_{k=0}^{p-1}\f{\b{2k}k^2\b{3k}k}{15^{3k}}\mod p.$$ Now applying
[S2, Theorem 4.6] we obtain the result.

\pro{Theorem 3.7} Let $p>3$ be a prime. Then
$$\align &\sum_{n=0}^{p-1}\f{D_n}{(-2)^n}
\e \sum_{n=0}^{p-1}\f{D_n}{4^n} \e
 \sum_{n=0}^{p-1}\f{D_n}{16^n}
 \e \sum_{n=0}^{p-1}\f{D_n}{(-32)^n}
\\&\e\sum_{k=0}^{p-1}\f{\b{2k}k^2\b{3k}k}{1458^k}
\e\cases 4x^2\mod p&\t{if $p\e 1\mod 3$ and so $p=x^2+3y^2$,}
\\0\mod p&\t{if $p\e 2\mod 3$.}
\endcases\endalign$$
\endpro
Proof. Taking $u=-\f 12,\f 1{16}$ in Theorem 3.1 and $u=-\f 12,\f
14,-\f 1{32}$ in Theorem 3.5 we see that
$$\align \sum_{n=0}^{p-1}\f{D_n}{(-2)^n}
\e \sum_{n=0}^{p-1}\f{D_n}{4^n} \e
 \sum_{n=0}^{p-1}\f{D_n}{16^n}
 \e \sum_{n=0}^{p-1}\f{D_n}{(-32)^n}
\e\sum_{k=0}^{p-1}\f{\b{2k}k^2\b{3k}k}{1458^k}
\e\sum_{k=0}^{p-1}\f{\b{2k}k^2\b{3k}k}{108^k}\mod p.\endalign$$ From
[M] and [Su2] we know that
$$\sum_{k=0}^{p-1}\f{\b{2k}k^2\b{3k}k}{108^k}
\e \cases 4x^2-2p\mod {p^2}&\t{if $p\e 1\mod 3$ and so
$p=x^2+3y^2$,}
\\0\mod {p^2}&\t{if $p\e 2\mod 3$.}
\endcases$$
Thus the result follows.
\par{\bf Remark 3.4} Let $p>5$ be a prime. In [S1] the author
conjectured that
$$\sum_{k=0}^{p-1}\f{\b{2k}k^2\b{3k}k}{1458^k}
\e \cases 4x^2-2p\mod {p^2}&\t{if $p\e 1\mod 3$ and so
$p=x^2+3y^2$,}
\\0\mod {p^2}&\t{if $p\e 2\mod 3$.}
\endcases$$
and
$$\sum_{k=0}^{p-1}\f{\b{2k}k^2\b{3k}k}{15^{3k}}
\e \cases 4x^2-2p\mod {p^2}&\t{if $p\e 1,4\mod {15}$ and so
$p=x^2+15y^2$,}\\2p-12x^2\mod {p^2}& \t{if $p\e 2,8\mod {15}$ and so
$p=3x^2+5y^2$,}
\\0\mod {p^2}&\t{if $p\e 7,11,13,14\mod {15}$.}
\endcases$$
In [Su3], Z.W. Sun conjectured that
$$\align&\sum_{n=0}^{p-1}\f {D_n}{(-2)^n}
\e \sum_{n=0}^{p-1}\f {D_n}{4^n}\e \sum_{n=0}^{p-1}\f {D_n}{16^n}
\e\sum_{n=0}^{p-1}\f {D_n}{(-32)^n}
\\&\e\cases 4x^2-2p\mod{p^2}&\t{if $p\e 1\mod 3$ and so $p=x^2+3y^2$,}
\\0\mod{p^2}&\t{if $p\e 2\mod 3$.}
\endcases\endalign$$

\section*{4. Congruences involving $\{a_n\}$}
\par For any nonnegative integer $n$ let
 $a_n=\sum_{k=0}^n\b nk^2\b{2k}k$.
Using Maple and the Zeilberger algorithm we find that
$$(n+2)^2a_{n+2}-(10n^2+30n+23)a_{n+1}+9(n+1)^2a_n=0\q
(n=0,1,2,\ldots).\tag 4.1$$
 \pro{Lemma 4.1} For any nonnegative
integer $n$ we have
$$\sum_{k=0}^n\b{2k}k\b{n+k}{2k}(-9)^{n-k}a_k
=\sum_{k=0}^n\b{2k}k^2\b{4k}{2k}\b{n+3k}{4k}(-27)^{n-k}.$$
\endpro
Proof. Let
$$\align &S_1(n)=\sum_{k=0}^n\b{2k}k\b{n+k}{2k}(-9)^{n-k}a_k,
\\&S_2(n)=\sum_{k=0}^n\b{2k}k^2\b{4k}{2k}\b{n+3k}{4k}(-27)^{n-k}.\endalign$$
Then $S_1(0)=1=S_2(0)$ and $S_1(1)=-3=S_2(1)$. Using the Maple
software doublesum.mpl and the method in [CHM] we find that for
$i=1,2$ and $m=0,1,2,\ldots$,
$$(m+2)^3S_i(m+2)+(2m+3)(7m^2+21m+17)S_i(m+1)+81(m+1)^3S_i(m)=0.$$
Thus $S_1(n)=S_2(n)$. This proves the lemma.

 \pro{Theorem 4.1} Let
$p$ be an odd prime and $u\in\Bbb Z_p$ with $(9u+1)(27u+1)\not\e
0\mod p$. Then
$$\sum_{k=0}^{p-1}\b{2k}k\Ls u{(1+9u)^2}^ka_k
\e\sum_{k=0}^{p-1}\b{2k}k^2\b{4k}{2k}\Ls u{(1+27u)^4}^k\mod p.$$
 \endpro
Proof. As $p\mid \b{2k}k$ for $\f p2<k<p$, we see that
$$\align &\sum_{k=0}^{p-1}\b{2k}k\Ls u{(1+9u)^2}^ka_k
\\&\e \sum_{k=0}^{(p-1)/2}\b{2k}ka_ku^k(1+9u)^{p-1-2k}
\e\sum_{k=0}^{(p-1)/2}\b{2k}ka_ku^k\sum_{r=0}^{p-1-2k}\b{p-1-2k}r(9u)^r
\\&\e\sum_{n=0}^{p-1}u^n\sum_{k=0}^n\b{2k}ka_k\b{p-1-2k}{n-k}9^{n-k}
\e\sum_{n=0}^{p-1}u^n\sum_{k=0}^n\b{2k}ka_k\b{-1-2k}{n-k}9^{n-k}
\\&=\sum_{n=0}^{p-1}u^n \sum_{k=0}^n\b{2k}ka_k(-9)^{n-k}
\b{n+k}{2k} \mod p.\endalign$$
On the other hand, as $p\mid
\b{2k}k\b{4k}{2k}$ for $\f p4<k<p$ we see that
$$\align &\sum_{k=0}^{p-1}\b{2k}k^2\b{4k}{2k}\Ls u{(1+27u)^4}^k
\\&\e \sum_{k=0}^{[p/4]}\b{2k}k^2\b{4k}{2k}u^k(1+27u)^{p-1-4k}
=\sum_{k=0}^{[p/4]}\b{2k}k^2\b{4k}{2k}u^k
\sum_{r=0}^{p-1-4k}\b{p-1-4k}r(27u)^r
\\&\e \sum_{n=0}^{p-1}u^n\sum_{k=0}^n\b{2k}k^2\b{4k}{2k}
\b{p-1-4k}{n-k}27^{n-k} \e
\sum_{n=0}^{p-1}u^n\sum_{k=0}^n\b{2k}k^2\b{4k}{2k}
\b{-1-4k}{n-k}27^{n-k}
\\&=\sum_{n=0}^{p-1}u^n
\sum_{k=0}^n\b{2k}k^2\b{4k}{2k}(-27)^{n-k} \b{n+3k}{4k} \mod
p.\endalign$$ Now combining all the above with Lemma 4.1 we deduce
the result.

\pro{Theorem 4.2} Let $p$ be a prime such that $p\e \pm 1\mod 8$.
Then
$$\sum_{k=0}^{p-1}\f{\b{2k}ka_k}{36^k}
\e \cases 4x^2\mod p&\t{if $p\e 1,7\mod{24}$ and so $p=x^2+6y^2$,}
\\0\mod p&\t{if $p\e 17,23\mod{24}$.}
\endcases$$
\endpro
Proof. Taking $u=\f 19$ in Theorem 4.1 and then applying [S3,
Theorem 5.4] we obtain the result.

\pro{Theorem 4.3} Let $p$ be a prime such that $p\e \pm 1\mod {12}$.
Then
$$\sum_{k=0}^{p-1}\f{\b{2k}ka_k}{(-12)^k}
\e \cases 4x^2\mod p&\t{if $p\e 1\mod{12}$ and so $p=x^2+9y^2$,}
\\0\mod p&\t{if $p\e 11\mod{12}$.}
\endcases$$
\endpro
Proof. Taking $u=-\f 13$ in Theorem 4.1 and then applying [S3,
Theorem 5.3] we obtain the result.

\section*{5. Congruences involving $\{b_n\}$}
\par Let
$\{b_n\}$ be the Almkvist-Zudilin numbers  given by (1.3).
  Since $\b mk\b kr=\b mr\b{m-r}{k-r}$, we see that
 $$ \align &\b{3k}k\b n{3k}\b{n+k}k\\&=\b nk\b{n-k}{2k}\b{n+k}k
 =\b {n-k}{2k}\b {2k}k\b{n+k}{2k}
 =\b{2k}k\b{4k}{2k}\b{n+k}{4k}.\endalign$$
Thus, $$ b_n=\sum_{k=0}^{[n/3]}
\b{2k}k^2\b{4k}{2k}\b{n+k}{4k}(-3)^{n-3k}.\tag 5.1$$
 From [CZ, Corollary 4.3] we know that
$$b_n=\sum_{k=0}^n\b{2k}k^2\b{4k}{2k}\b{n+3k}{4k}(-27)^{n-k}.
\tag 5.2$$ This is true since  $b_n$ and
$b_n'=\sum_{k=0}^n\b{2k}k^2\b{4k}{2k}\b{n+3k}{4k}(-27)^{n-k}$ have
the same initial values and recurrence relation:
$$(n+2)^3b_{n+2}+(2n+3)(7n^2+21n+17)b_{n+1}+81(n+1)^3b_n=0.$$

\pro{Theorem 5.1} Let $p$ be an odd prime and $u\in\Bbb Z_p$ with
$u\not\e -\f 1{27}\mod p$. Then
$$\sum_{n=0}^{p-1}b_nu^n
\e \sum_{k=0}^{p-1}\b{2k}k^2\b{4k}{2k}\Ls u{(1+27u)^4}^k\mod p.$$
\endpro
Proof. As $p\mid \b{2k}k\b{4k}{2k}$ for $\f p4<k<p$, using Fermat's
little theorem and (5.2) we see that
$$\align \sum_{k=0}^{p-1}\b{2k}k^2\b{4k}{2k}\Ls u{(1+27u)^4}^k
&\e \sum_{k=0}^{[p/4]}\b{2k}k^2\b{4k}{2k}u^k(1+27u)^{p-1-4k}
\\&=\sum_{k=0}^{[p/4]}\b{2k}k^2\b{4k}{2k}u^k\sum_{r=0}^{p-1-4k}
\b{p-1-4k}r(27u)^r
\\&\e\sum_{n=0}^{p-1}u^n\sum_{k=0}^n\b{2k}k^2\b{4k}{2k}\b{p-1-4k}{n-k}
 27^{n-k}
 \\&\e \sum_{n=0}^{p-1}u^n\sum_{k=0}^n\b{2k}k^2\b{4k}{2k}
 \b{-1-4k}{n-k}27^{n-k}
 \\&=\sum_{n=0}^{p-1}u^n\sum_{k=0}^n\b{2k}k^2\b{4k}{2k}
 \b{n+3k}{4k}(-27)^{n-k}
 \\&=\sum_{n=0}^{p-1}b_nu^n\mod p.\endalign$$
 This proves the theorem.
 \pro{Corollary 5.1} Let
$p$ be an odd prime and $u\in\Bbb Z_p$ with $(9u+1)(27u+1)\not\e
0\mod p$. Then
$$\sum_{k=0}^{p-1}\b{2k}k\Ls u{(1+9u)^2}^ka_k\e \sum_{n=0}^{p-1}
b_nu^n\mod p.$$
\endpro
Proof. This is immediate from Theorems 4.1 and 5.1.

\pro{Theorem 5.2} Let $p>3$ be a prime. Then
 $$\sum_{n=0}^{p-1}\f{b_n}{(-9)^n}\e
 \cases 4x^2\mod p&\t{if $p\e 1\mod 3$ and so $p=x^2+3y^2$,}
 \\0\mod p&\t{if $p\e 2\mod 3$.}
 \endcases$$
 \endpro
 Proof. Taking $u=-\f 19$ in Theorem 5.1 and then
  applying [S3, Theorem 5.1] we deduce the result.

 \pro{Theorem 5.3} Let $p$ be a prime with $p\e \pm 1\mod 8$. Then
 $$\sum_{n=0}^{p-1}\f{b_n}{9^n}\e
 \cases 4x^2\mod p&\t{if $p\e 1,7\mod {24}$ and so $p=x^2+6y^2$,}
 \\0\mod p&\t{if $p\e 17,23\mod {24}$.}
 \endcases$$
 \endpro
 Proof. Taking $u=\f 19$ in Theorem 5.1 and then
 applying [S3, Theorem 5.4] we deduce the result.

\pro{Theorem 5.4} Let $p$ be an odd prime and $u\in\Bbb Z_p$ with
$u\not\e -\f 13\mod p$. Then
$$\sum_{n=0}^{p-1}b_nu^n
\e \sum_{k=0}^{p-1}\b{2k}k^2\b{4k}{2k}\Ls {u^3}{(1+3u)^4}^k\mod p.$$
\endpro
Proof. As $p\mid \b{2k}k\b{4k}{2k}$ for $\f p4<k<p$, using Fermat's
little theorem and (5.1) we see that
$$\align \sum_{k=0}^{p-1}\b{2k}k^2\b{4k}{2k}\Ls {u^3}{(1+3u)^4}^k
&\e \sum_{k=0}^{[p/4]}\b{2k}k^2\b{4k}{2k}u^{3k}(1+3u)^{p-1-4k}
\\&=\sum_{k=0}^{[p/4]}\b{2k}k^2\b{4k}{2k}u^{3k}\sum_{r=0}^{p-1-4k}
\b{p-1-4k}r(3u)^r
\\&\e\sum_{n=0}^{p-1}u^n\sum_{k=0}^{[n/3]}\b{2k}k^2\b{4k}{2k}
\b{p-1-4k}{n-3k}
 3^{n-3k}
 \\&\e \sum_{n=0}^{p-1}u^n\sum_{k=0}^{[n/3]}\b{2k}k^2\b{4k}{2k}
 \b{-1-4k}{n-3k}(-3)^{n-3k}
 \\&=\sum_{n=0}^{p-1}u^n\sum_{k=0}^{[n/3]}\b{2k}k^2\b{4k}{2k}
 \b{n+k}{4k}(-3)^{n-3k}
 \\&=\sum_{n=0}^{p-1}b_nu^n\mod p.\endalign$$
 This proves the theorem.

 \pro{Corollary 5.2} Let $p$ be an odd prime and $u\in\Bbb Z_p$ with
$u\not\e -\f 13,-\f 1{27}\mod p$. Then
$$\sum_{k=0}^{p-1}\b{2k}k^2\b{4k}{2k}\Ls u{(1+27u)^4}^k
\e \sum_{k=0}^{p-1}\b{2k}k^2\b{4k}{2k}\Ls {u^3}{(1+3u)^4}^k\mod p.$$
\endpro
Proof. This is immediate from Theorems 5.1 and 5.4.

 \pro{Theorem 5.5} Let $p$ be a prime with $p\e \pm
1\mod{12}$. Then
 $$\sum_{n=0}^{p-1}\f{b_n}{(-3)^n}\e \sum_{n=0}^{p-1}\f{b_n}{(-27)^n}
 \e\cases 4x^2\mod p&\t{if $p\e 1\mod {12}$ and so $p=x^2+9y^2$,}
 \\0\mod p&\t{if $p\e 11\mod {12}$.}
 \endcases$$
 \endpro
 Proof. Taking $u=-\f 13$ in Theorem 5.1 and $u=-\f 1{27}$
 in Theorem 5.4
 we see that
 $$\sum_{n=0}^{p-1}\f{b_n}{(-3)^n}\e \sum_{k=0}^{p-1}
 \f{\b{2k}k^2\b{4k}{2k}}{(-12288)^k}
 \e \sum_{n=0}^{p-1}\f{b_n}{(-27)^n}\mod p.$$
 Now applying [S3, Theorem 5.3] we deduce the result.

\pro{Theorem 5.6} Let $p$ be a prime such that $p>7$. Then
$$\sum_{k=0}^{p-1}\f{\b{2k}k^2\b{4k}{2k}}{28^{4k}}
\e\cases 4x^2\mod p&\t{if $p=x^2+2y^2\e 1,3\mod 8$,}
 \\0\mod {p^2}&\t{if $p\e 5,7\mod 8$}
 \endcases$$
 and
 $$\align \sum_{k=0}^{p-1}\f{\b{2k}ka_k}{100^k}
 \e \sum_{n=0}^{p-1}b_n\e
 \sum_{n=0}^{p-1}\f{b_n}{81^n}
 \e\cases 4x^2\mod p&\t{if $p=x^2+2y^2\e 1,3\mod 8$,}
 \\0\mod p&\t{if $p\e 5,7\mod 8$.}
 \endcases\endalign$$
 \endpro
 Proof. Putting $u=1$ in Theorems 4.1 and 5.1
 and $u=\f 1{81}$ in Theorem 5.4
 we see that
 $$ \sum_{k=0}^{p-1}\f{\b{2k}ka_k}{100^k}
 \e\sum_{n=0}^{p-1}b_n\e \sum_{n=0}^{p-1}\f{b_n}{81^n}
 \e  \sum_{k=0}^{p-1}\f{\b{2k}k^2\b{4k}{2k}}{28^{4k}}\mod p.$$
 Taking $u=1$ in Corollary 5.2 we see that
 $$\sum_{k=0}^{p-1}\f{\b{2k}k^2\b{4k}{2k}}{28^{4k}}\e
 \sum_{k=0}^{p-1}\f{\b{2k}k^2\b{4k}{2k}}{256^k}\mod p.$$
From [M] and [Su2] we know that
$$\sum_{k=0}^{p-1}\f{\b{2k}k^2\b{4k}{2k}}{256^k}
\e\cases 4x^2-2p\mod {p^2} &\t{if $p=x^2+2y^2\e 1,3\mod 8$,}
\\0\mod {p^2}&\t{if $p\e 5,7\mod 8$.}
\endcases$$
Thus, the result is true when the modulus is $p$. By [S3, Theorem
4.2],
$$\sum_{k=0}^{p-1}\f{\b{2k}k^2\b{4k}{2k}}{28^{4k}} \e
P_{[\f p4]}\Ls{20\sqrt 6}{49}^2\mod p.$$ Hence, using [S3, Theorem
4.2] again we see that for primes $p\e 5,7\mod 8$,
$$P_{[\f p4]}\Ls{20\sqrt 6}{49}\e 0\mod p\qtq{and so}
\sum_{k=0}^{p-1}\f{\b{2k}k^2\b{4k}{2k}}{28^{4k}} \e 0\mod{p^2}.$$
The proof is now complete.
\par {\bf Remark 5.1} In [Su1], Zhi-Wei Sun conjectured that for any prime
$p\not=2,3, 7$,
$$\sum_{k=0}^{p-1}\f{\b{2k}k^2\b{4k}{2k}}{28^{4k}}
\e \cases 4x^2-2p\mod {p^2}&\t{if $p=x^2+2y^2\e 1,3\mod 8$,}
\\0\mod {p^2}&\t{if $p\e 5,7\mod 8$.}
\endcases$$

 \section*{6. Some conjectures
on congruences modulo prime powers}
\par For any nonnegative integers $n$ let $\{D_n\},\{a_n\}$ and
$\{b_n\}$ be given by
 (1.3).   Suppose that $p>3$ is a prime. In [Su4] Z.W. Sun conjectured congruences for
$\sum_{k=0}^{p-1} \f{a_k}{3^k}$ and $\sum_{k=0}^{p-1}
\f{a_k}{(-3)^k}\mod {p^2}$. In [Su3, Conjecture 7.8] Z.W. Sun
conjectured explicit congruences for $\sum_{k=0}^{p-1}
\b{2k}k\f{a_k}{m^k}\mod {p^2}$ in the cases $m=198^2,-1123596$.
 \par By doing calculations with the help of Maple,
 we pose some conjectures. These conjectures are similar to
 some conjectures in [Su1, Su3, S1]. As showed in [S1-S4], many
 conjectures for supercongruences are connected with binary
 quadratic forms of class number $1$ or $2$ and the number of points
 on certain elliptic curves with complex multiplication
 over the field $\Bbb F_p$ with $p$
 elements.

\pro{Conjecture 6.1} Let $p$ be a prime with $p\e
1,17,19,23\mod{30}$. Then
$$\sum_{k=0}^{p-1}D_k
\e\sum_{k=0}^{p-1}\f {D_k}{64^k}
\e\sum_{k=0}^{p-1}\f{\b{2k}k^2\b{3k}k}{(-27)^k}\mod{p^3}.$$
\endpro

\pro{Conjecture 6.2} Let $p$ be a prime greater than $3$. Then
$$\sum_{k=0}^{p-1}\f {D_k}{(-8)^k}
\e \sum_{k=0}^{p-1}\f{\b{2k}k^2\b{3k}k}{216^k}\mod{p^3} \qtq{for}p\e
1,5,7,11\mod{24}$$ and
$$\sum_{k=0}^{p-1}\f {D_k}{8^k}
\e \sum_{k=0}^{p-1}\f{\b{2k}k^2\b{3k}k}{8^k}\mod{p^3}\qtq{for}p\e
1,3\mod 8.$$
\endpro

\pro{Conjecture 6.3} Let $p>3$ be a prime. Then
$$\sum_{k=0}^{p-1}\f {D_k}{4^k}
\e\sum_{k=0}^{p-1}\f {D_k}{(-32)^k} \e \Ls
p3\sum_{k=0}^{p-1}\f{\b{2k}k^2\b{3k}k}{108^k}\mod{p^3},$$ and for
$p\e 1\mod 3$ we have
$$\align\sum_{k=0}^{p-1}\f {D_k}{(-2)^k}
\e \sum_{k=0}^{p-1}\f {D_k}{16^k} \e
\sum_{k=0}^{p-1}\f{\b{2k}k^2\b{3k}k}{108^k} \e
\sum_{k=0}^{p-1}\f{\b{2k}k^2\b{3k}k}{1458^k} \e
\sum_{k=0}^{p-1}\f{\b{2k}k^3}{16^k} \mod{p^3}.\endalign$$
\endpro

\pro{Conjecture 6.4} Let $p>3$ be a prime. Then
$$\sum_{k=0}^{p-1}\f{\b{2k}ka_k}{36^k}\e
\sum_{k=0}^{p-1}\f{b_k}{9^k} \e \cases 4x^2-2p\mod {p^2}&\t{if
$p=x^2+6y^2\e 1,7\mod{24}$,}
\\ 2p-8x^2\mod{p^2}&\t{if $p=2x^2+3y^2\e 5,11\mod{24}$,}
\\0\mod {p^2}&\t{if $p\e 13,17,19,23\mod{24}$.}
\endcases$$
\endpro
\pro{Conjecture 6.5} Let $p>5$ be a prime. Then
$$\sum_{k=0}^{p-1}\f{\b{2k}ka_k}{100^k}\e
\sum_{k=0}^{p-1}b_k\e \sum_{k=0}^{p-1}\f{b_k}{81^k}  \e \cases
4x^2-2p\mod {p^2}&\t{if $p=x^2+2y^2\e 1,3\mod 8$,}
\\0\mod {p^2}&\t{if $p\e 5,7\mod 8$.}
\endcases$$
\endpro

\pro{Conjecture 6.6} Let $p>3$ be a prime. Then
$$\align&\sum_{k=0}^{p-1}\f{\b{2k}ka_k}{(-12)^k}
\e \sum_{k=0}^{p-1}\f{b_k}{(-3)^k}\e
\sum_{k=0}^{p-1}\f{b_k}{(-27)^k}
\\&\e \cases 4x^2-2p\mod {p^2}&\t{if
$12\mid p-1$ and so $p=x^2+9y^2$,}
\\2p-2x^2\mod{p^2}&\t{if $12\mid p-5$ and so $2p=x^2+9y^2$,}
\\0\mod {p^2}&\t{if $p\e 3\mod 4$.}
\endcases\endalign$$
\endpro

\pro{Conjecture 6.7} Let $p>3$ be a prime. Then
$$\sum_{k=0}^{p-1}\f{b_k}{(-9)^k}
\e\cases 4x^2-2p\mod{p^2}&\t{if $p=x^2+3y^2\e 1\mod 3$,}
\\0\mod{p^2}&\t{if $p\e 2\mod 3$.}
\endcases$$\endpro

\pro{Conjecture 6.8} Let $p>3$ be a prime. Then
$$\align&\sum_{k=0}^{p-1}(4k+1)\f{b_k}{(-27)^k}
\e \sum_{k=0}^{p-1}(4k+1)\f{b_k}{81^k}\e
\sum_{k=0}^{p-1}(2k+1)\f{b_k}{(-9)^k} \e
\sum_{k=0}^{p-1}(2k+1)\f{b_k}{9^k}
\\&\e \f 13\sum_{k=0}^{p-1}(4k+3)b_k
\e \f 13\sum_{k=0}^{p-1}(4k+3)\f{b_k}{(-3)^k} \e \Ls p3p\mod{p^2}.
\endalign$$
\endpro

\pro{Conjecture 6.9} For any prime $p\e 5\mod 6$ we have
$\sum_{k=0}^{p-1}\f{\b{2k}k^2\b{3k}k}{1458^k}\e 0\mod{p^3}.$
\endpro
\pro{Conjecture 6.10} Let $p>3$ be a prime. Then
$$\sum_{n=0}^{p-1}\sum_{k=0}^n\b nk^4\e \Ls p3\sum_{k=0}^{p-1}
\f{\b{2k}k^2\b{3k}k}{(-27)^k}\mod{p^3}.$$
\endpro
\par Z.W. Sun made a conjecture on $\sum_{n=0}^{p-1}
\sum_{k=0}^n\b nk^4\mod{p^2}$. See Remark 3.1.

 \pro{Conjecture 6.11} Let $p>5$ be a prime.
Then
$$\align &\sum_{k=0}^{p-1}\f{63k+8}{(-15)^{3k}}
\b{2k}k\b{3k}k\b{6k}{3k}\e 8p\Ls {-15}p\mod{p^3},
\\&\sum_{k=0}^{p-1}\f{133k+8}{255^{3k}}
\b{2k}k\b{3k}k\b{6k}{3k}\e 8p\Ls {-255}p\mod{p^3}\qtq{for}p\not=17,
\\&\sum_{k=0}^{p-1}\f{28k+3}{20^{3k}}
\b{2k}k\b{3k}k\b{6k}{3k}\e 3p\Ls{-5}p\mod{p^3},
\\&\sum_{k=0}^{p-1}\f{63k+5}{66^{3k}}
\b{2k}k\b{3k}k\b{6k}{3k}\e 5p\Ls{-33}p\mod{p^3}\qtq{for}p\not=11,
\\&\sum_{k=0}^{p-1}\f{11k+1}{54000^k}
\b{2k}k\b{3k}k\b{6k}{3k}\e p\Ls {-15}p\mod{p^3},
\\&\sum_{k=0}^{p-1}\f{506k+31}{(-12288000)^k}
\b{2k}k\b{3k}k\b{6k}{3k}\e 31p\Ls{-30}p\mod{p^3}.
\endalign$$
\endpro
\par Conjecture 6.11 is similar to some conjectures in [Su1].

\end{document}